\documentclass[11pt]{article}
  \usepackage[active]{srcltx} 
  \usepackage{amsmath}
  \usepackage{amscd}
  \usepackage{amsfonts}
  \usepackage{amssymb}
  \usepackage{amsthm}
  \usepackage{euscript}
  \usepackage{graphicx}
  \usepackage{epsfig}
  \usepackage{dsfont}
  \usepackage{lscape}

\theoremstyle{definition}

\begin{document}
\begin{center}
{\Large\bf 
The Centroid Method for Imaging through Turbulence
\par\vspace{.2cm}
 (Preliminary version)
}
\par\vspace{.7cm}
{
Mario Micheli
}
\par\vspace{.25cm}
\normalsize
MAP5, Universit\'e Paris Descartes
\\ 
45, rue des Saints P\`eres, 7\`eme \'etage 
\\
75270 Paris Cedex 06, France
\\ \tt mariomicheli@gmail.com 
%
\end{center}
\par
\vspace*{.2cm}
\begin{abstract}
A simple and effective 
method for imaging through ground-level atmospheric turbulence.
\end{abstract}
\section{Introduction}
The problem of reconstructing an image from a sequence that is altered by
ground-level atmospheric turbulence is still largely unsolved,
due to the time- and space-dependent 
nature of the distortion. In fact what is commonly 
referred to as ÒturbulenceÓ in imaging can be modeled, at least in first approximation, as the combined effect of (i) a time-dependent deformation of 
the image domain and (ii) a blur with an anisoplanatic point spread function.
This suggests an approach to solve the problem by
first correcting for the geometric distortion and later applying a deblurring algorithm.
\par
Such approach is followed by other authors; see, for example,~\cite{gilles,gilles:cam1,shimizu:08,milanfar:10,milanfar:11}. In several cases, however, the geometric correction is achieved by simply computing the temporal mean of the data sequence. This introduces 
additional blur, whose effects, in the case
of strong turbulence, are particularly destructive of image features. 
In this report we present a novel \em geometric \em
idea, which is the one of correcting for the domain distortion by  
first computing the temporal mean of the \em deformations \em of the images with respect to the first
one, and then creating the \em centroid \em image by warping
the first one via this average deformation. At a second stage the data (i.e.~the image sequence)
is registered onto the centroid and then averaged in time; this produces a blurri image that is made sharp by nonlocal Total Variation (NL-TV) deconvolution~\cite{gilboa:2,yifei:10}. 
In order to formulate all of this mathematically 
we need a notion of deformation; with this in mind, we briefly summarize the concept of {\em optical flow\/} and the standard techniques for its recovery.

\section{Optical flow}
Roughly speaking, \em optical flow \em 
between two images is 
``a vector field that deforms one image into the other''. More precisely, 
given an image sequence~$I(x,y,t)$, $(x,y)\in\Omega, t\in[0,T]$, 
the optical flow 
between the images $I(x,y,t_0)$
and
$I(x,y,t_0+\delta t)$, $(x,y)\in\Omega$,
(for fixed~$t_0$ and~$\delta t$) 
is a vector field $\mathbf{ f}(x,y)=(u(x,y),v(x,y))$ 
such that:
$I\big(x+u(x,y),y+v(x,y),t_0+\delta t\big) \simeq I(x,y,t_0)$, $(x,y)\in \Omega$.
\par
There is a very large and consolidated literature on the recovery of optical flow; see, for example, \cite{barron,black,horn:1,lucaskanade,singh,tomasikanade} and references therein. A good part of it relies on the \em data correlation 
constraint\em\/, i.e.~on the assumption that
the image brightness of a moving region 
essentially remains constant in time. 
That is,
if the location of a point  moves within the image frame according to  
the map~\mbox{$
t\mapsto(x(t),y(t))$}, 
then it is the case that~$\frac{d}{dt}I(x(t),y(t),t)=0$.
The chain rule yields:
$
\frac{\partial I}{\partial x}
\frac{dx}{dt}
+
\frac{\partial I}{\partial y}
\frac{dy}{dt}
+
\frac{\partial I}{\partial t}=0,
$
where the velocity field~$(u,v):=(dx/dt,dy/dt)$
is interpreted as the optical flow at time~$t$ and location~$(x,y)$.
If we rewrite such equation as
\begin{equation}
\label{transport}
u
\frac{\partial I}{\partial x}
+
v
\frac{\partial I}{\partial y}
+
\frac{\partial I}{\partial t}=0,
\end{equation}
which must be solved in~$u$ and~$v$, we 
note that the recovery of optical flow
may be regarded as an \em inverse problem\em~\cite{aster, kirsch}. In fact while
equation~\eqref{transport} is the \em transport equation\em~\cite{evans} for the function~$I$,
it is 
parameters~$(u,v)$ that must be recovered (the function~$I$ is known from the data). Like most inverse problems
optical flow recovery is \em ill-posed \em
in that at each location $(x,y)\in\Omega$
the solution of~\eqref{transport} is a line in the~$(u,v)$-plane.
Moreover, in the texture free regions (characterized by $\partial I/\partial x=\partial I/\partial y=0$)
the problem is even more underdetermined.
Therefore further hypotheses must be introduced.
\par
A fast and effective method was introduced by 
Horn~\& Schunk~\cite{horn:1,horn:2}:
given an image sequence~$I$
the optical flow \em at a fixed time~$t$ \em
is a vector field~$(u,v):\Omega\rightarrow \mathbb{R}^2$
that minimizes the energy
\begin{equation}
\label{HSe}
E[u,v]
=\iint_\Omega
\Big\{
\Big(
u\frac{\partial I}{\partial x} 
+
v\frac{\partial I}{\partial y}
+
\frac{\partial I}{\partial t}
\Big)^2
+
\alpha
\|\nabla u \|^2
+
\alpha
\|\nabla v \|^2
\Big\}\,dx\,dy;
\end{equation}
in other words we have added the $H^1$-seminorm
of the optical flow~$(u,v)$ as regularizing term.
The Euler-Lagrange equations of~\eqref{HSe} are the partial
differential equations:
\begin{align*}
\frac{\partial I}{\partial x}
\Big(
\frac{\partial I}{\partial x} u
+
\frac{\partial I}{\partial y} v
+
\frac{\partial I}{\partial t}
\Big)
-\alpha\Delta u
&=
0,
&
\frac{\partial I}{\partial y}
\Big(
\frac{\partial I}{\partial x} u
+
\frac{\partial I}{\partial y} v
+
\frac{\partial I}{\partial t}
\Big)
-\alpha\Delta v
&=
0.
\end{align*}
One sees immediately that in a texture-free region,
in which the first terms on the left-hand side of the 
above equations are zero,
the resulting optic flow is a vector field with
{\it harmonic} components, whose boundary values are determined by 
the image regions that do have texture. Note also that
the parameter~$\alpha>0$ essentially determines the  smoothness of the solution~$(u,v)$.
\par
Other optical flow recovery methods have emerged: most notably,
Black \& Anandan~\cite{black,black:2010} 
have developed a robust scheme whose output
is sharp and allows one to detect the boundaries of moving objects.
In our application the motion field is produced by atmospheric
turbulence and is thus relatively smooth by its own nature: whence  
Horn \& Schunk's algorithm, which has the further advantage  
of being computationally efficient,
will be our method of choice.
\section{The Centroid}
\label{sec_centroid}
\begin{figure}[t]
\begin{center}
\hspace*{-3cm}
\begin{picture}(300,80)
\setlength{\unitlength}{1.5pt}
\put(10,2){\makebox(0,0){$P_1$}}
\put(14,46){\makebox(0,0){$P_2$}}
\put(58,52){\makebox(0,0){$P_3$}}
\put(78,30){\makebox(0,0){$P_4$}}
\put(39,36){\makebox(0,0){$c$}}
\put(10,10){\makebox(0,0){$\bullet$}}
\put(20,40){\makebox(0,0){$\bullet$}}
\put(50,50){\makebox(0,0){$\bullet$}}
\put(70,30){\makebox(0,0){$\bullet$}}
\put(40,30){\makebox(0,0){$\circ$}}
\put(110,2){\makebox(0,0){$P_1$}}
\put(114,46){\makebox(0,0){$P_2$}}
\put(158,52){\makebox(0,0){$P_3$}}
\put(178,30){\makebox(0,0){$P_4$}}
\put(144,34){\makebox(0,0){ $c$}}
\put(110,10){\makebox(0,0){$\bullet$}}
\put(120,40){\makebox(0,0){$\bullet$}}
\put(150,50){\makebox(0,0){$\bullet$}}
\put(170,30){\makebox(0,0){$\bullet$}}
\put(140,30){\makebox(0,0){$\circ$}}
\put(210,2){\makebox(0,0){$P_1$}}
\put(214,46){\makebox(0,0){$P_2$}}
\put(258,52){\makebox(0,0){$P_3$}}
\put(278,30){\makebox(0,0){$P_4$}}
\put(244,28){\makebox(0,0){$c$}}
\put(238,45){\makebox(0,0){ $q$}}
\put(210,10){\makebox(0,0){$\bullet$}}
\put(220,40){\makebox(0,0){$\bullet$}}
\put(250,50){\makebox(0,0){$\bullet$}}
\put(270,30){\makebox(0,0){$\bullet$}}
\put(240,40){\makebox(0,0){$\bullet$}}
\put(240,30){\makebox(0,0){$\circ$}}
\thinlines
\qbezier[20](10,10)(15,25)(20,40)
\qbezier[20](20,40)(35,45)(50,50)
\qbezier[18](70,30)(60,40)(50,50)
\qbezier[40](10,10)(40,20)(70,30)
%
\qbezier[20](120,40)(135,45)(150,50)
\qbezier[18](170,30)(160,40)(150,50)
%
\qbezier[20](210,10)(215,25)(220,40)
\qbezier[20](220,40)(235,45)(250,50)
\qbezier[18](270,30)(260,40)(250,50)
\qbezier[40](210,10)(240,20)(270,30)
\put(10,10){\line(3,2){29}} %
\put(20,40){\line(2,-1){18.8}} %
\put(50,50){\line(-1,-2){9.25}} %
\put(70,30){\line(-1,0){28.8}} %
\put(110,10){\vector(1,3){9.57}} %
\put(110,10){\vector(1,1){39}} %
\put(110,10){\vector(3,1){58.6}} %
%
\put(240,40){\vector(1,1){9}} %
\put(240,40){\vector(3,-1){29}} %
\put(240,40){\vector(-1,0){18.8}} %
\put(240,40){\vector(-1,-1){28.8}} %
\thicklines
\put(110,10.2){\vector(3,2){28.9}} %
\put(110.1,9.8){\vector(3,2){28.9}} %
\put(110,10){\vector(3,2){28.9}} %
\put(240,40){\vector(0,-1){8.6}} %
\put(239.75,40){\vector(0,-1){8.6}} %
\end{picture}
\\
\hspace*{.9cm}
(a)
\hspace{4.5cm}
(b)
\hspace{4.5cm}
(c)
\end{center}
\vspace*{-.5cm}
\caption{(a) The centroid in Euclidean space.
(b) $\protect\overrightarrow{P_1C}$ (thick)
is obtained by averaging vectors~$\protect\overrightarrow{P_1P_i}$ (thin).
(c) The ``correction'' $\protect\overrightarrow{QC}$ (thick) is obtained by averaging
vectors $\protect\overrightarrow{QP_i}$ (thin). 
}
\label{c_ill}
\end{figure}
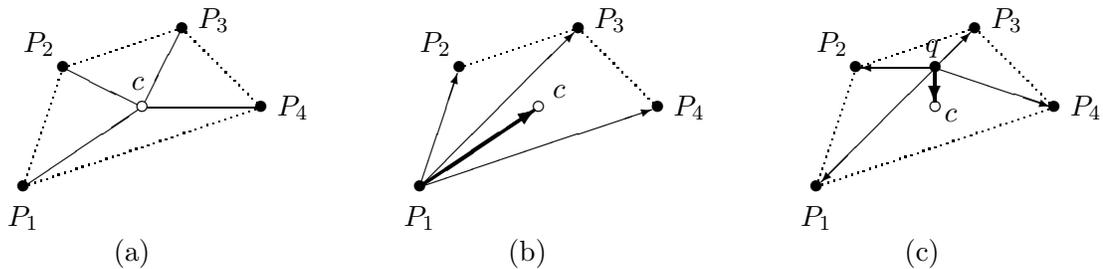

In Euclidean geometry the \em centroid \em 
of a set of points~$\{P_1,\ldots,P_N\}$
whose coordinates are known 
is given 
by~$C=\frac{1}{N}\sum_{i=1}^NP_i$; 
see Figure~\ref{c_ill}(a).
If only the positions of the points \em with respect to
one of them\em, e.g.~without loss of generality~$P_1$, are known, i.e.~if  only vectors
$\overrightarrow{P_1P_i}$, $i=2,\ldots N$ are known, then
we can still compute the position 
of~$C$ \em with respect to~$P_1$ \em
by calculating
\begin{equation}
\label{centr_eq}
\overrightarrow{P_1C}
=
\frac{1}{N}
\sum_{i\not=1}
\overrightarrow{P_1P_i}\,,
\end{equation}
see Figure~\ref{c_ill}(b).
Also, it is immediate to verify that
$
\sum_{i=1}^N \overrightarrow{CP_i}=0
$,
and that 
$
\overrightarrow{QC}=\frac{1}{N}\sum_{i=1}^N\overrightarrow{QP_i}
$
for any point~$Q$; see Figure~\ref{c_ill}(c).
\par
These elementary considerations lead to the following idea.
Suppose we want to compute the ``centroid'' of the set of
$N$ images~$\{I_i(x), x\in\Omega\}_{1\leq i\leq N}$. Thinking of 
images as
``points'' and of the optical flow between two images
as a ``vector'' between them, we can compute
the optical flow between the reference image~$I_1$ 
and the rest of them, 
i.e.~the $N$ ``vectors'' $\overrightarrow{I_1I_i}(x)$, 
$x\in\Omega$, $i=1,\ldots,N$;
so, treating optical flow as an element of a linear space, 
we compute
the average flow~$\overrightarrow{I_1C}$ 
using
equation~\eqref{centr_eq}.
This guiding principle may be summarized as follows:
$$
\boxed{
\begin{array}{c}
\mbox{\it Instead of computing the temporal mean of the images,}
\\ 
\mbox{\em we compute the temporal mean of the {\bf deformations}}
\\ 
\mbox{\em of the images with respect
to the first one.}
\end{array}
}
$$
\par
Once this average deformation (optical flow)~$\overrightarrow{I_1C}$ 
is known, and~we warp 
$I_1$ via the flow~$\overrightarrow{I_1C}$ to finally get
the {\it centroid image}~$C$ of the data set.
Note however that since we are dealing with images (and not 
points in Euclidean space) 
this procedure is not exact.
Most notably,
the resulting~$C$ actually depends quite strongly on the choice of the initial reference image: if certain features are missing in~$I_1$
(due to the effects of optical turbulence), they will also
be absent in the centroid, since the latter is obtained by warping~$I_1$.
Also, it is the case that  
 $\sum_{i=1}^N\overrightarrow{C I_i}$ is usually \em not \em equal to zero.
In fact, a ``better'' result is obtained by implementing
the two following procedures.
\par 
\paragraph{\bf A. Iteration.}
The following iterative algorithm
can be implemented:
\begin{enumerate}
\item a first estimate of the centroid is computed with the procedure described above, and we shall call this~$Q$;
\item the \em correction \em  flow
$
\overrightarrow{QC}:=\frac{1}{N}\sum_{i=1}^N\overrightarrow{QP_i}
$
is calculated (typically this is a non-zero vector field), see Figure~\ref{c_ill}(c); a new estimate of the centroid image~$C$
is computed by warping~$Q$ via the optical 
flow~$\overrightarrow{QC}$;
\item the previous step is repeated until convergence is achieved.
\end{enumerate}
In fact the correction will never become zero, but its~$L^2$ norm will decrease at each step. Our experience with our data sets is that such norm typically stabilizes after about 4 or 5 iterations. 
\paragraph{\bf B. Registration of the data onto the Centroid, and averaging.} The image that is obtained by the previous procedure
is better than the initial estimate of the centroid (i.e.~the one 
obtained with only one iteration), in that the geometric correction
is improved. However, as we noted above, since the centroid
is obtained by warping one particular image of the data set,~say~$I_1$,
it will lack the features that the optical turbulence has 
eliminated in such frame. A procedure to alleviate this problem is the following: 
\begin{enumerate}
\item \em register \em each image~$I_i$ onto the centroid~$C$; 
this can be done by computing 
the flow~$\varphi_i:=\overrightarrow{CI_i}$
and defining~$I_i^\mathrm{R}(x):=I_i(x+\varphi(x))$, $x\in\Omega$.
\item \em average \em the resulting registered sequence,
$
C^\mathrm{R}(x):=\frac{1}{N}\sum_{i=1}^NI_i^\mathrm{R}(x)$, $x\in\Omega$.
\end{enumerate}
Under the assumption that all features of interest appear in the majority of the data frames, the above procedure allows one to 
recover them (albeit in blurred form) in the 
resulting image~$C^\mathrm{R}$.
\section{Deblurring}
So far we have used the data sequence~$I_1,\ldots,I_N$
to produce an image~$C^\mathrm{R}$ where the geometric
distortion has been rectified, but that still is blurry.
In other words, we have reduced our problem to a \em deblurring problem\em\/, for which we employ known deconvolution techniques.
\par
Given the fact that in most of our images there are repeated patterns,
our method of choice is
\em nonlocal Total Variation \em deconvolution
(NL-TV); see~\cite{gilboa:2,yifei:10}. The energy functional is of the type
\begin{equation}
\label{NLTVe}
E[u]
:=\int_{\Omega}(f - k \ast u)^2 dx
+\alpha\int_\Omega
\sqrt{
\int_\Omega
\big[
u(x)-u(y)
\big]^2
w_f(x,y)
\,dy}
\,dx
\end{equation}
where the weight is defined as 
$$
w_f(x,y)
:=
\exp\Big(
-\frac{\big(G_a\ast|f(x+\cdot)-f(y+\cdot)|\big)(0)}{h^2}
\Big).
$$
with $x,y\in\Omega$. In the above formula~$G_a$
is a Gaussian mask with standard deviation~$a$ and $h$
is a scaling parameter. 
Note that the weight is approximately equal to~1 if the
patches around locations~$x$ and~$y$ are similar,
otherwise it is approximately zero. By the 
second term on the right-hand side of~\eqref{NLTVe}, that
the weight
``enforces'' the same similarities that appear in the data~$f$
onto the minimizers~$u$ of the above energy.
\par
In our case the data~$f$ is the blurry image~$C^\mathrm{R}$, whereas the~$k$
is a Guassian kernel whose standard deviation is chosen manually in a way that yields 
the most visually appealing result. In the future we plan to test physics-inspired 
point-spread functions, such as the Fried kernel described in~\cite{gilles:cam2}.  
\section{Results}
\begin{figure}
\begin{center}
\mbox{ }\par\vspace{.2cm}
\begin{picture}(0,0)
\put(-3,60){\makebox(0,0)[r]{(A)}}
\end{picture}
\includegraphics[height=4.10cm]{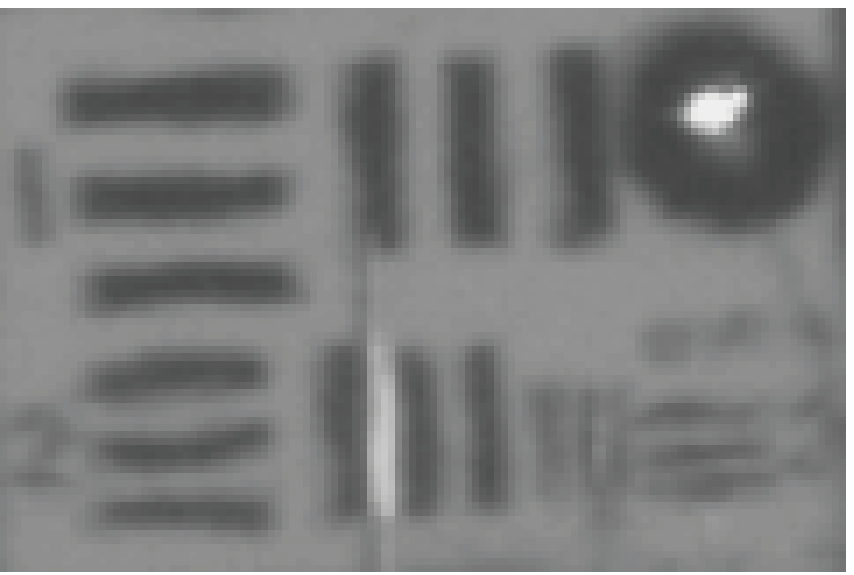}%
\includegraphics[height=4.10cm]{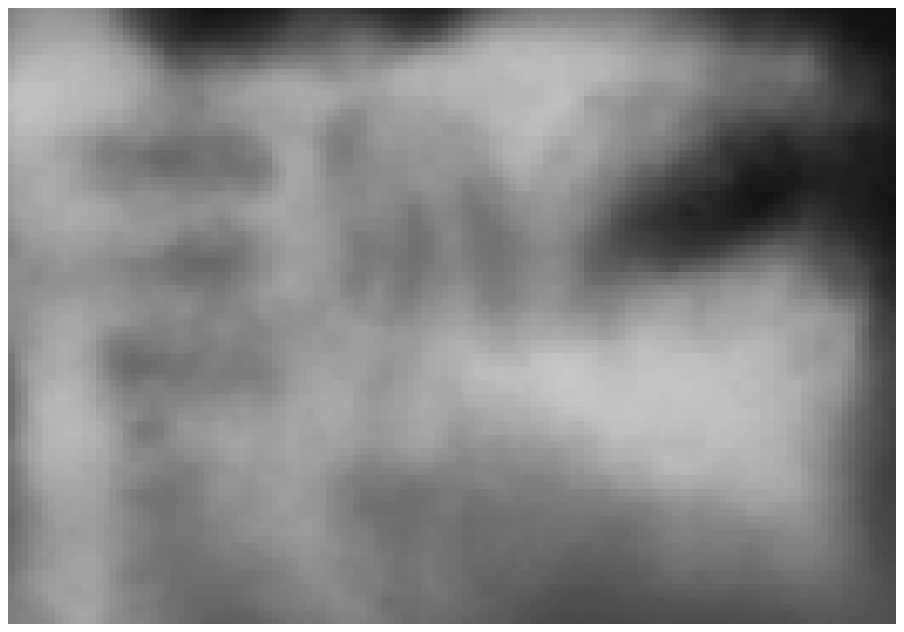}%
\includegraphics[height=4.10cm]{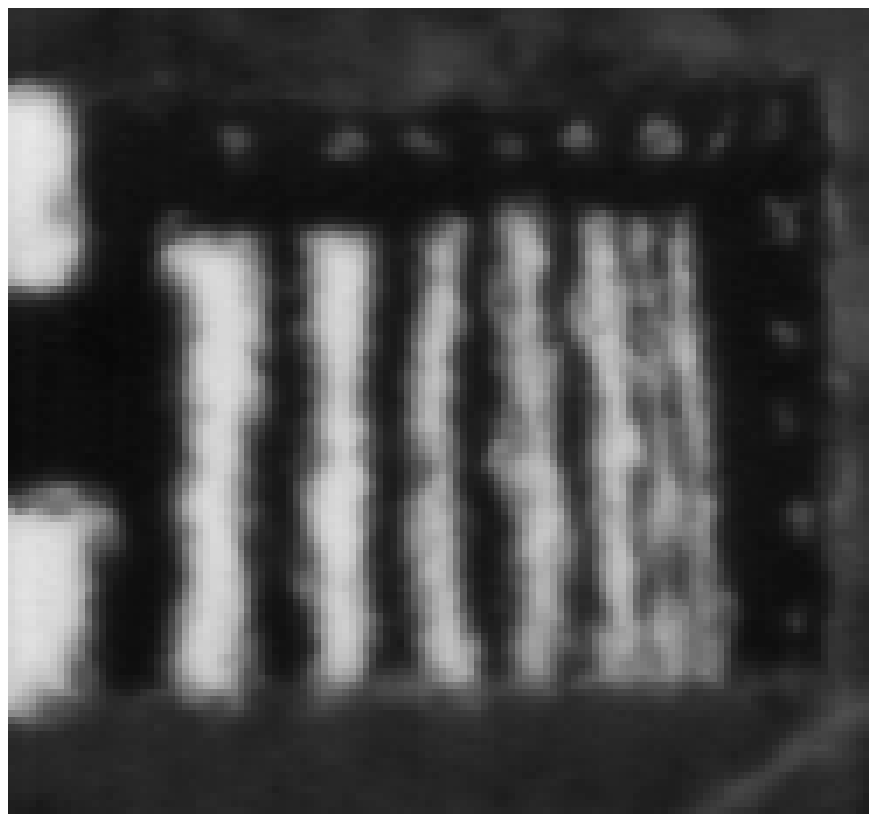}
\\
\begin{picture}(0,0)
\put(-3,60){\makebox(0,0)[r]{(B)}}
\end{picture}
\includegraphics[height=4.10cm]{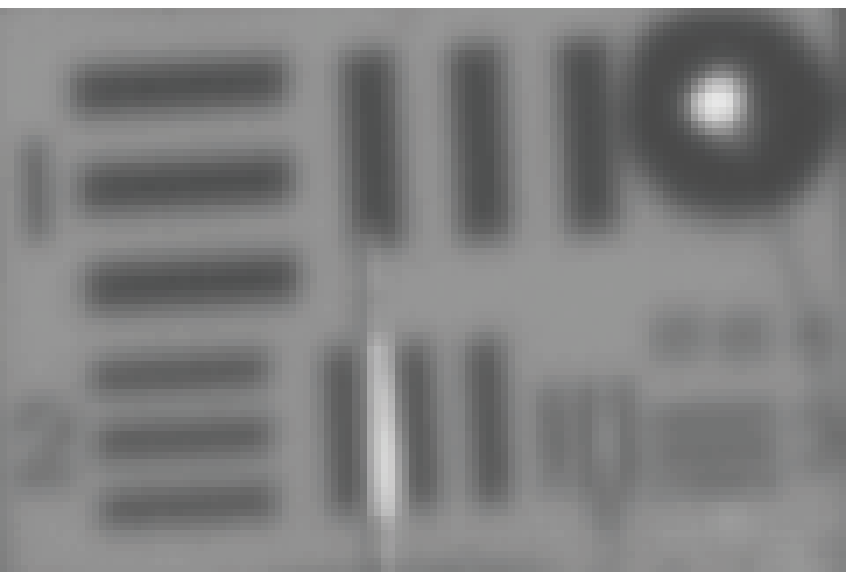}%
\includegraphics[height=4.10cm]{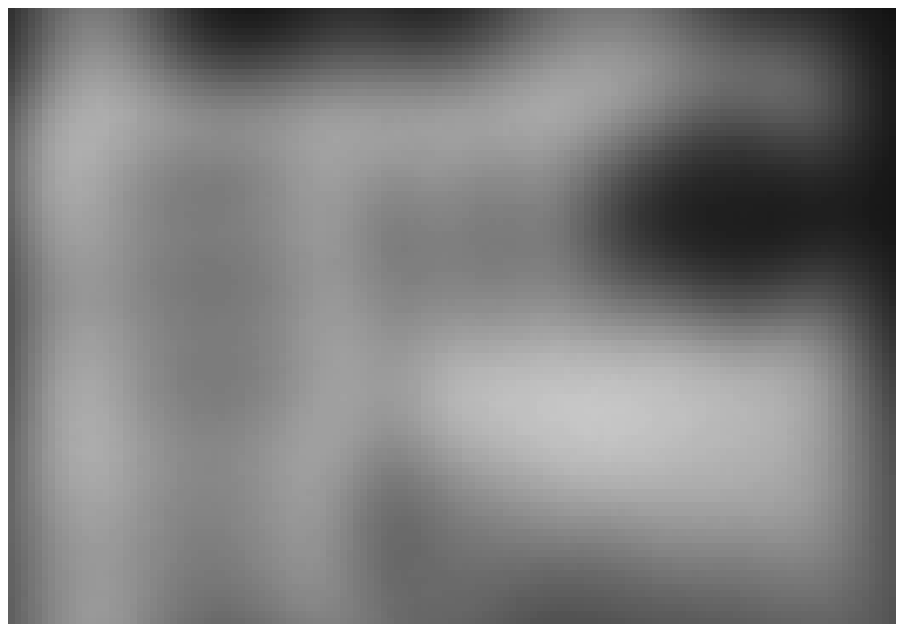}%
\includegraphics[height=4.10cm]{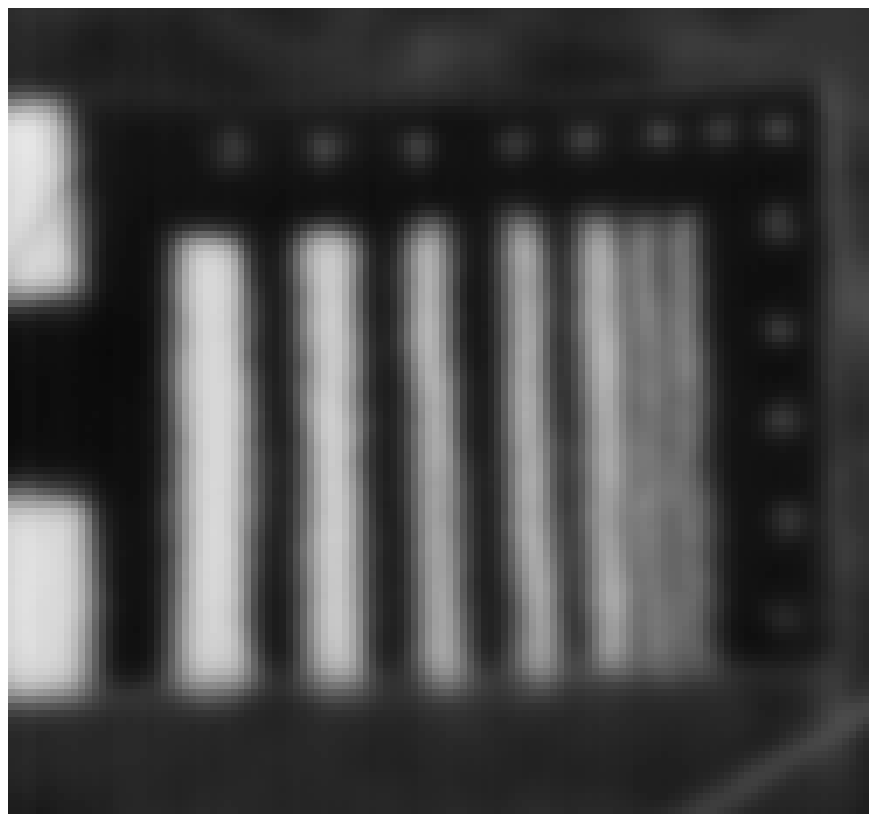}
\\
\begin{picture}(0,0)
\put(-3,60){\makebox(0,0)[r]{(C)}}
\end{picture}
\includegraphics[height=4.10cm]{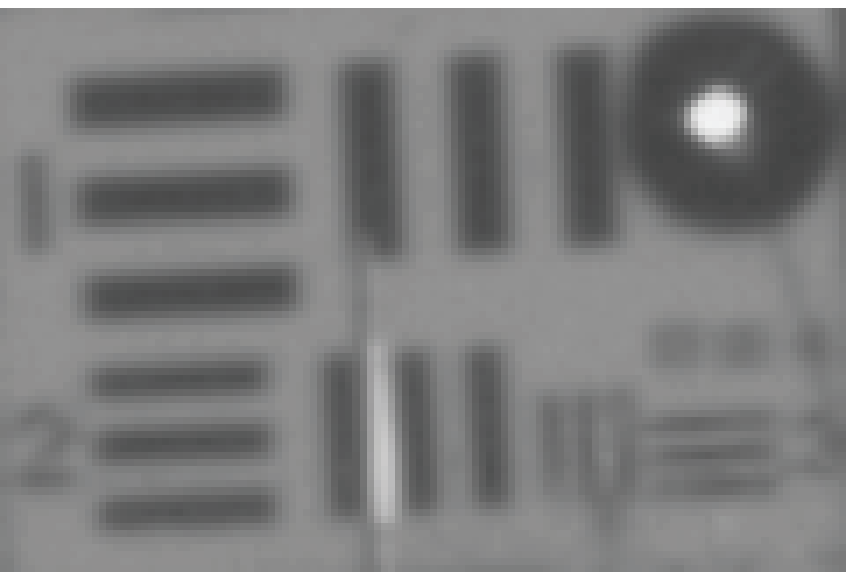}%
\includegraphics[height=4.10cm]{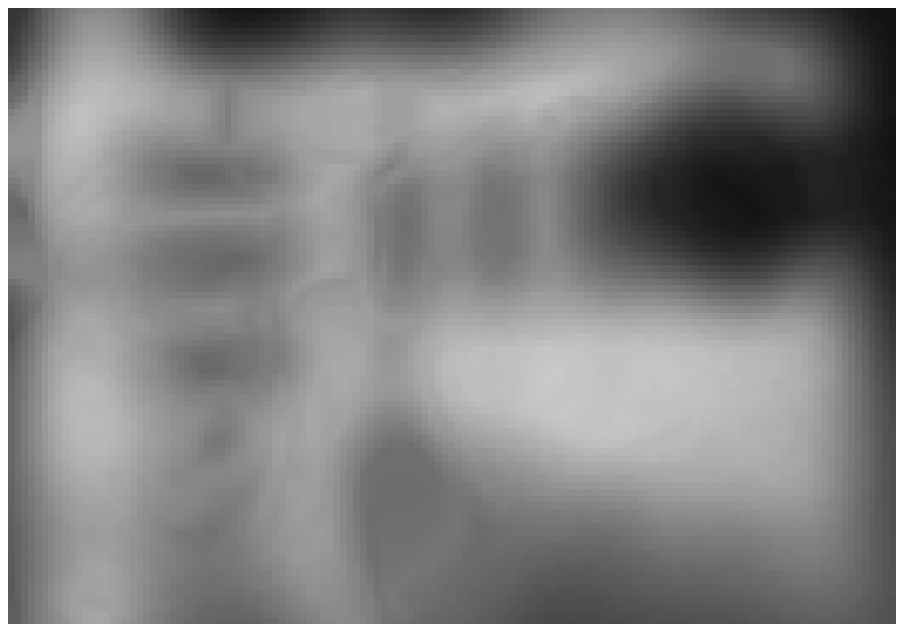}%
\includegraphics[height=4.10cm]{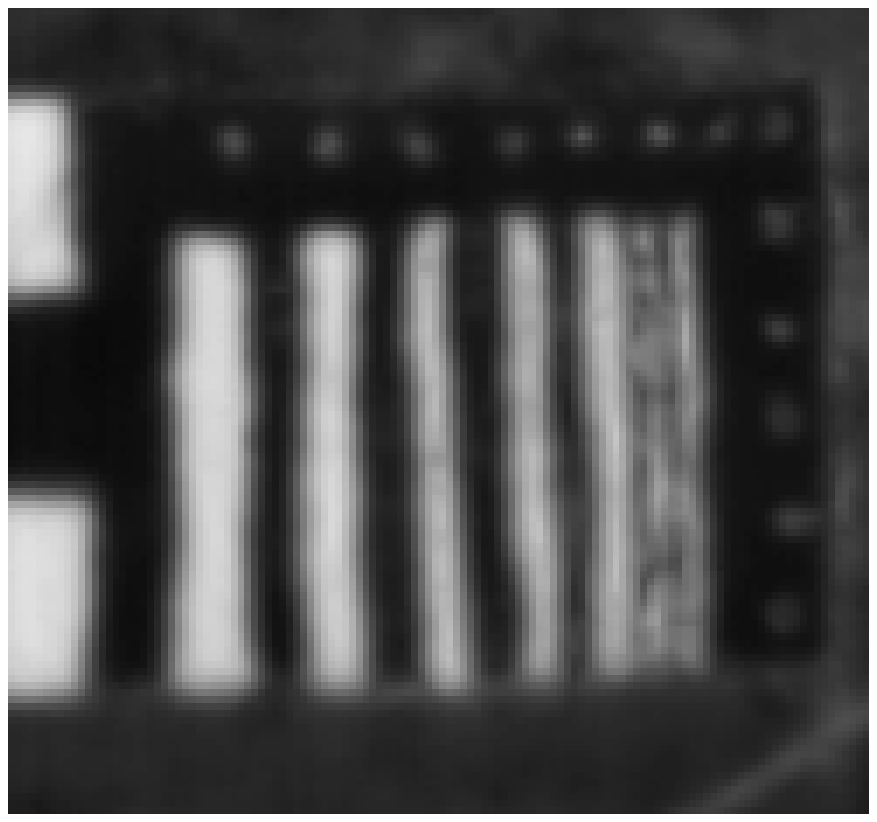}
\\
\begin{picture}(0,0)
\put(-3,60){\makebox(0,0)[r]{(D)}}
\end{picture}
\includegraphics[height=4.10cm]{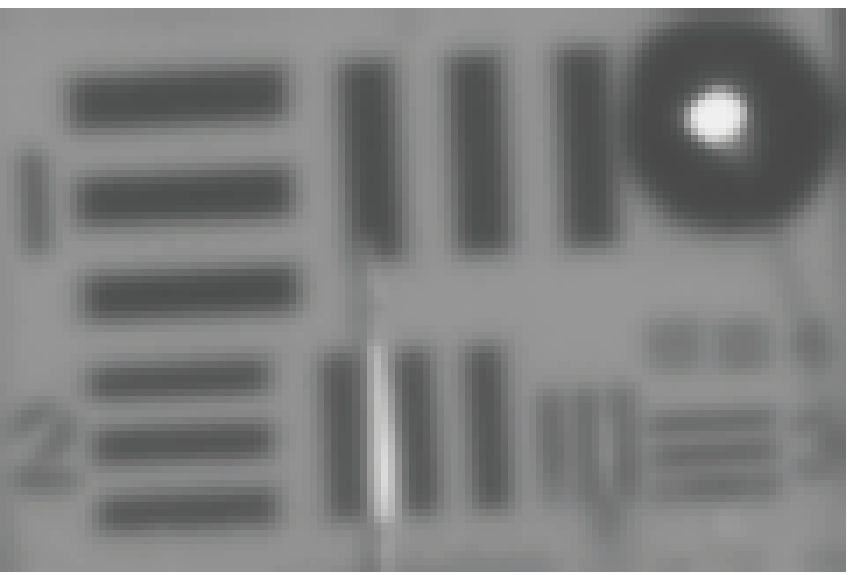}%
\includegraphics[height=4.10cm]{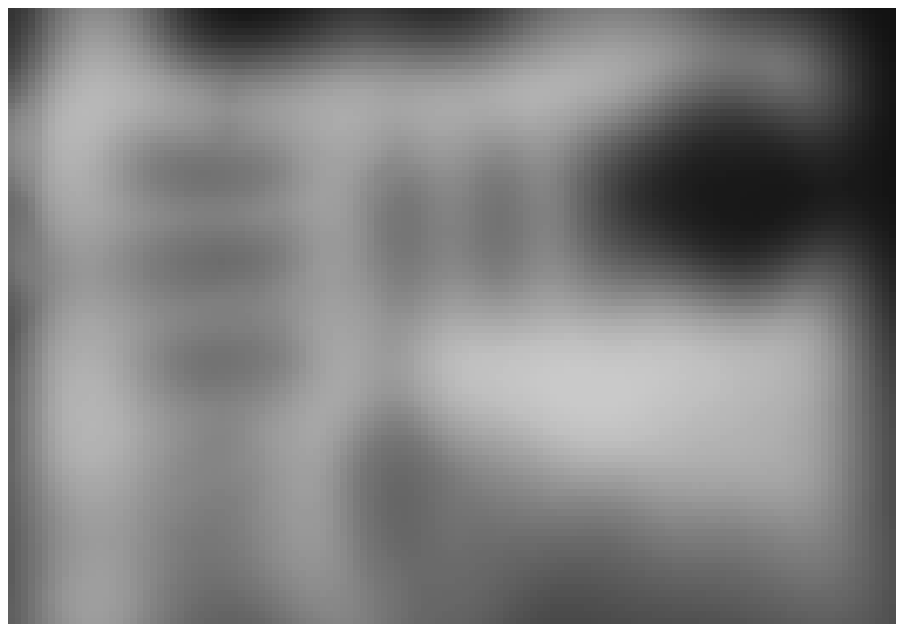}%
\includegraphics[height=4.10cm]{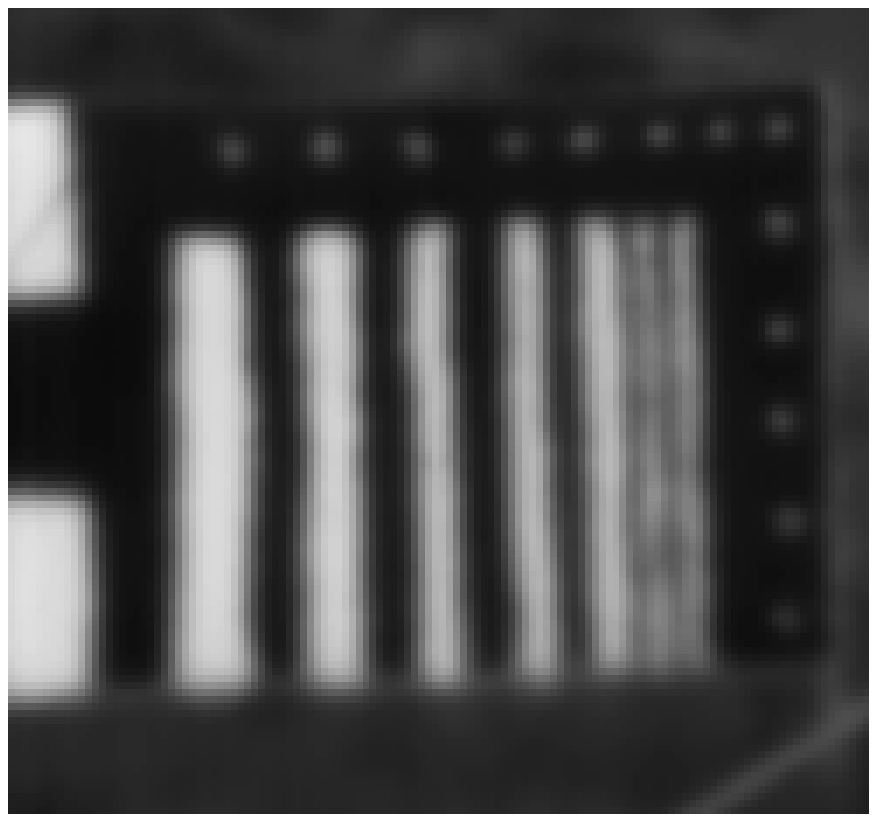}
\\
\begin{picture}(0,0)
\put(-3,60){\makebox(0,0)[r]{(E)}}
\end{picture}
\includegraphics[height=4.10cm]{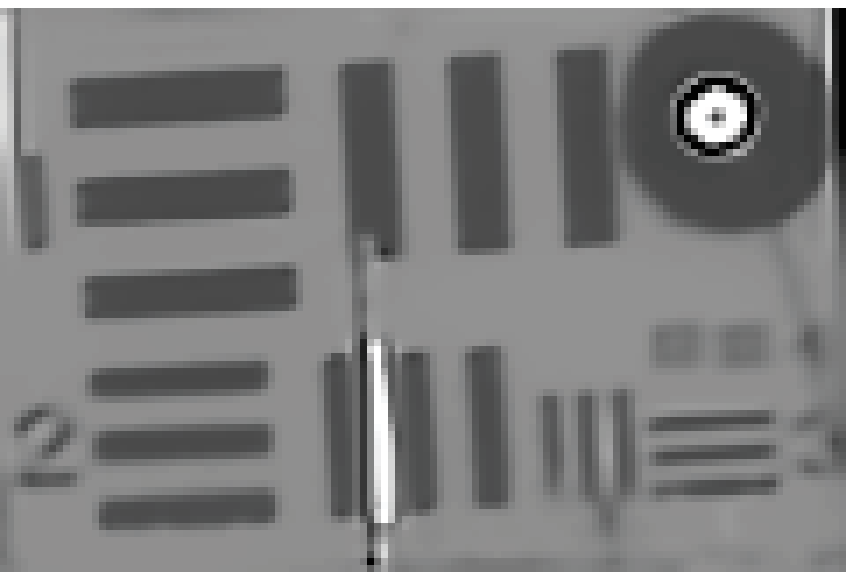}%
\includegraphics[height=4.10cm]{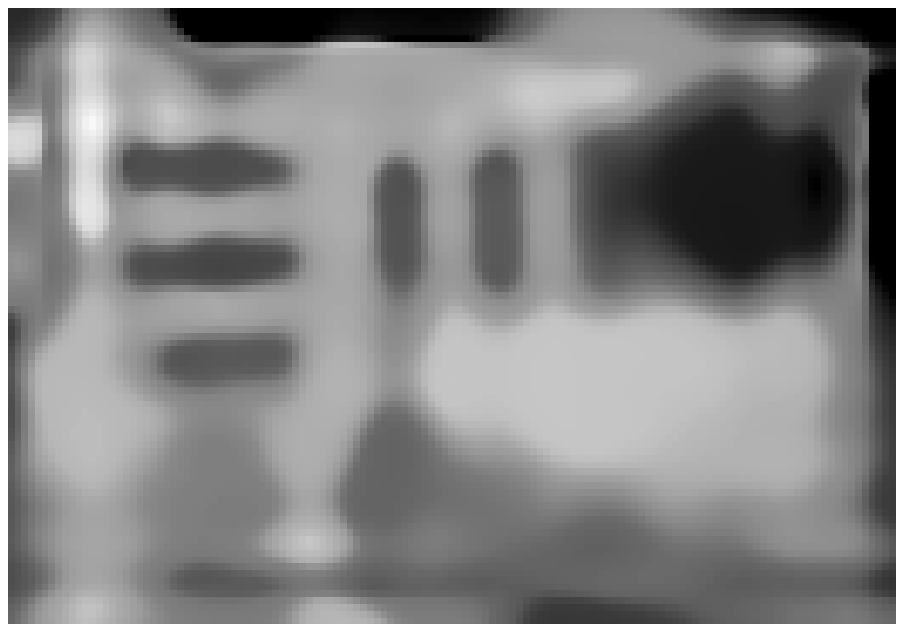}%
\includegraphics[height=4.10cm]{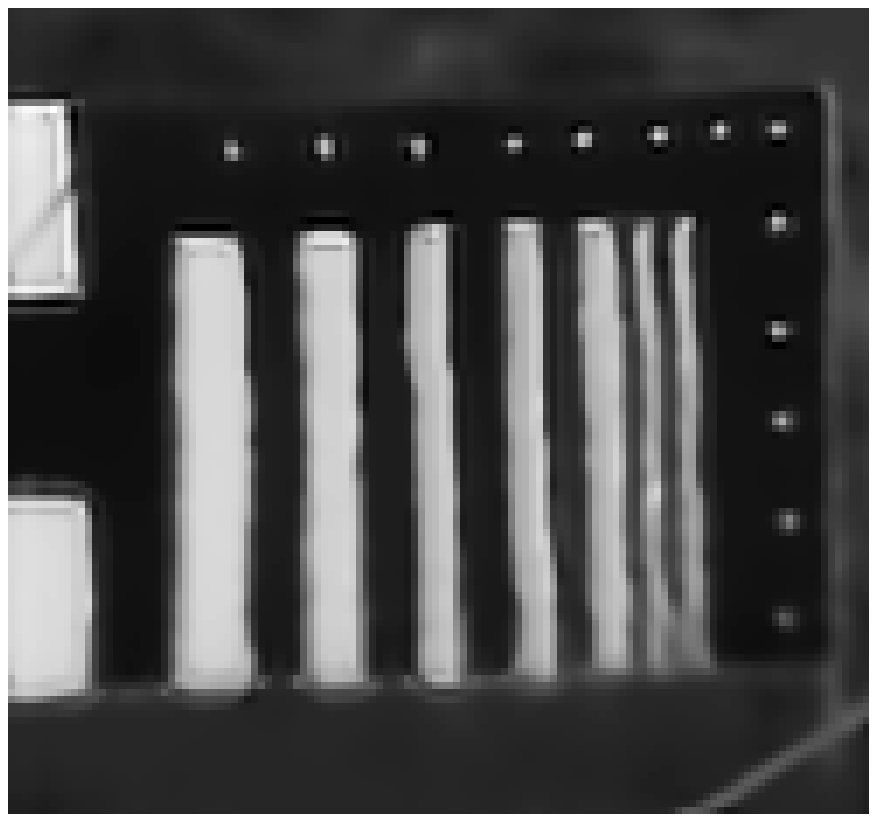}
\end{center}
\par\vspace{-.3cm}
\hspace{2.7cm}(1)\hspace{5.6cm}(2)\hspace{4.6cm}(3)
\caption{Results for three image sequences.}
\label{fig_res}
\end{figure}

In Figure~\ref{fig_res} we present the results for three data sets.
For each one we show: 
\par\vspace{.2cm}
(A) a frame from the original data sequence; 
\par
(B) the temporal mean;
\par
(C) the centroid~$C$;
\par
(D) the mean~$C^\mathrm{R}$ of the images 
obtained by registering the data set onto the centroid;
\par
(E) the image resulting from applying NL-TV to~$C^\mathrm{R}$.
\par\vspace{.2cm}
In fact we never use the temporal mean in our algorithm, but we show it here to 
emphasize how the centroid is generally much sharper and
preserves more features than the temporal mean; this is especially true for data sets with
intense turbulence, such as the one shown in the second column of Figure~\ref{fig_res}.  
We should note that while the centroid contains artifacts that are
inherited by the particular image that is used to create it, this in not the case
with the image~$C^\mathrm{R}$ obtained by averaging the registered data.
This is particularly apparent for the NATO sequence, shown in
column~(3) in Figure~\ref{fig_res}: note in particular how the 
second bar from the right has much more regular edges in image (D)-(3)
than in image (C)-(3), and how
some of the dots are better defined.
\section{Current and future work}
The geometric method that we have briefly described here proves quite 
effective for many
data sets, such as the three that we have shown in this report. It is conceptually 
very simple and rather straightforward to implement.  The results are comparable
with the state of the art.
\par
Perhaps the most fundamental issue that is currently being explored is the dependence
of the centroid on the initial (reference) image, as we already discussed above; the procedure
described at the end of section~\ref{sec_centroid} (denoted with~{\bf B}) certainly 
alleviates the problem as it allows some of the possibly missing 
features to be recovered, but 
experimental results show that there is still a certain degree of dependence 
on the reference image. We are currently exploring algorithms to deal with this issue. 
Future work will also focus on speeding up the algorithm and selecting the parameters 
(e.g.~in the computation of optical flow, or in the deconvolution step) in an automated manner.\section{Acknowledgements}
The author's research was partly supported by 
ONR grant N000140910256
and the KaraMetria program of the Agence 
Nationale de la Recherche (ANR) of France.
I would like to thank Dr.~Alan Van Nevel at the U.S.~Naval 
Air Warfare Center, Weapons Division (China Lake, California) 
and the NATO SET156 (ex-SET072) Task Group
for providing the image data.
I am very grateful to  
Andrea Bertozzi,  J\'er\^ome Gilles, Yifei Lou, Stanley Osher, and Stefano Soatto
of UCLA for introducing me to this research area in the first place and for 
the many lengthy discussions on the topic. 
I would also like to express my gratitude to Joan Alexis Glaun\`es of Universit\'e Paris Descartes
for his advice and support during my stay at the MAP5 lab, and to
Gabriel Peyr\'e and Fran\c{c}ois-Xavier Vialard of Universit\'e Paris Dauphine for
their insightful suggestions. 
Last, but not least,  the author is indebted to Xiaoqun Zhang, formerly at UCLA and now at Shanghai Jiaotong University, for providing very efficient
nonlocal Total Variation (NL-TV) deconvolution code.  

\end{document}